\documentclass[12pt]{article}

\newcommand{\be}{\begin{equation}}
\newcommand{\ee}{\end{equation}}
\newcommand{\bi}[1]{\vspace{-3mm} \bibitem{#1}}

\usepackage{epsfig,amsmath,amssymb,graphics,graphicx} 

\begin{document}

\begin{center}

International Journal of Applied and Computational Mathematics. \\
2016. Vol.2. No.2. P.195-201. DOI: 10.1007/s40819-015-0054-6
\vskip 7mm

{\bf \large Local Fractional Derivatives of 
\vskip 3mm
Differentiable Functions are Integer-order Derivatives or Zero} \\

\vskip 7mm
{\bf \large Vasily E. Tarasov} \\
\vskip 3mm

{\it Skobeltsyn Institute of Nuclear Physics,\\ 
Lomonosov Moscow State University, Moscow 119991, Russia} \\
{E-mail: tarasov@theory.sinp.msu.ru} \\

\begin{abstract}
In this paper we prove that 
local fractional derivatives of differentiable functions
are integer-order derivative or zero operator.
We demonstrate that the local fractional derivatives
are limits of the left-sided Caputo fractional derivatives.
The Caputo derivative of fractional order $\alpha$
of function $f(x)$ is defined 
as a fractional integration of order $n-\alpha$
of the derivative $f^{(n)}(x)$ of integer order $n$. 
The requirement of the existence of integer-order derivatives
allows us to conclude that 
the local fractional derivative cannot be 
considered as the best method to describe 
nowhere differentiable functions and fractal objects.
We also prove that unviolated Leibniz rule 
cannot hold for derivatives of orders $\alpha \ne 1$.
\end{abstract}

\end{center}


\section{Introduction}

Fractional calculus as a theory of
derivatives and integrals of  non-integer orders
\cite{SKM,KST,Podlubny}
has a long history \cite{His1,His2,His3}. 
As a mathematical tool it has a lot of applications 
in physics and mechanics to
to describe media and systems with 
nonlocal and hereditary properties of 
of power-law type.
During more than three hundred years of history of fractional calculus it has been proposed different definitions
of fractional derivatives by Riemann, Liouville, Riesz, 
Gr\"unwald, Letnikov, Marchaud, Caputo \cite{SKM,KST}.
Fractional-order derivatives can be characterized
by a set of interesting unusual properties
such as a violation of the usual Leibniz rule, 
a deformations of the usual chain rule,
a violation of the semi-group property.
These unusual properties can be used  
to describe nontrivial and unusual properties of 
complex systems and media.
Therefore these properties are valuable and 
important for applications in physics and mechanics 
(see \cite{Mainardi,TarasovSpringer,Uch2} and 
references therein).


Recently, so-called local fractional derivative 
has been proposed in \cite{KG1996,KG1997,KG1998,Review}
by a motivation of describing
the properties of fractal objects and processes. 
The definition of the local fractional derivative 
is obtained from the Riemann-Liouville
fractional derivative by introducing two changes 
such as a subtraction of finite number of terms of 
the Taylor series and the taking the limit. 
The motivation of this definition is 
a studying fractional differentiability properties of 
nowhere differentiable functions and fractal objects. 

In this paper, we demonstrate that the definition
of local fractional derivative can be presented as
a limit of the Caputo fractional derivatives.
We note that the Caputo fractional derivative of 
order $\alpha$ ($n-1 <\alpha \le n$) of a function
is defined as an action of derivative of integer order $n$
and then the integration of non-integer order.
Therefore the function should be differentiable $n$ times.
In this paper we also prove that 
local fractional derivatives of differentiable functions
are integer-order derivative or zero operator.
We argued that unviolated Leibniz rule 
cannot hold for local fractional 
derivatives of orders $\alpha \ne 1$.


\section{Local fractional derivative via Caputo fractional derivative}

The local fractional derivative ${\cal D}^{\alpha}$ 
of order $\alpha$ of a function $f(x)$ 
at the point $x=a$ is usually defined 
(see Eq. 26 of \cite{KG1997})
by the equation
\be \label{LFD-n}
({\cal D}^{\alpha} f)(a) := \lim_{x \to a}
\, _a^{RL}{\cal D}^{\alpha}_x  
\left( f(x) - \sum^{n-1}_{k=0} \frac{f^{(k)}(a)}{k!} (x-a)^k \right) , 
\quad ( n-1< \alpha \le n)
\ee
where $ \, _a^{RL}{\cal D}^{\alpha}_x $ is
the left-sided Riemann-Liouville fractional derivative
of order $\alpha$ (see Eqs. (2.1.15) and (2.1.1) of \cite{KST})
that is defined by
\be \label{RLFD-n}
\, _a^{RL}{\cal D}^{\alpha}_x f(x) :=
\frac{1}{\Gamma(n-\alpha)} \frac{d^n}{dx^n} 
\int^x_a \frac{f(z) \, dz}{(x-z)^{\alpha-n+1}} .
\quad ( n-1< \alpha \le n)
\ee
The Riemann-Liouville fractional derivative
$ \, _a^{RL}{\cal D}^{\alpha}_x$
is nonlocal operator.
For $ 0< \alpha \le 1$, 
equations (\ref{LFD-n}) and (\ref{RLFD-n}) have the forms
\be \label{LFD-1}
({\cal D}^{\alpha} f)(a) := \lim_{x \to a}
\, _a^{RL}{\cal D}^{\alpha}_x  (f(x) - f(a)) , 
\quad ( 0< \alpha \le 1) ,
\ee
\be \label{RLFD-1}
\, _a^{RL}{\cal D}^{\alpha}_x f(x) :=
\frac{1}{\Gamma(1-\alpha)} \frac{d}{dx} 
\int^x_a \frac{f(z) \, dz}{(x-z)^{\alpha}} ,
\quad ( 0< \alpha \le 1) .
\ee

In the fractional calculus is well-known 
the Caputo fractional derivatives \cite{KST,Podlubny}
suggested in \cite{C1,C2,CM1,CM2}.
The left-sided Caputo fractional derivative
of order $\alpha$ can be defined 
via the above Riemann-Liouville fractional  derivative 
(see Eq. (2.4.1) and (2.4.3) of \cite{KST}) by
\be \label{CFD-n}
(\, _a^C{\cal D}^{\alpha}_x f)(x) := 
\, _a^{RL}{\cal D}^{\alpha}_x  
\left( f(x) - \sum^{n-1}_{k=0} \frac{f^{(k)}(a)}{k!} (x-a)^k \right) , 
\quad ( n-1< \alpha \le n) .
\ee
As a result, the local fractional derivative (\ref{LFD-n}) 
of order $\alpha$ at the point $x=a$ 
can be defined via the above Caputo fractional derivative 
by the equation
\be \label{LFD-nC}
({\cal D}^{\alpha} f)(a):= \lim_{x \to a}
\, (\, _a^C{\cal D}^{\alpha}_x f)(x) , 
\quad ( n-1< \alpha \le n) .
\ee
As a result, we can see that local fractional derivative 
is a limit of the Caputo fractional derivative. \\

{\bf Proposition.} \\
{\it The local fractional derivative (\ref{LFD-n}) 
of order $\alpha$, where $n-1 < \alpha \le n$ 
and $n \in \mathbb{N}$, 
at the point $x=x_0$ can be represented by the 
equation 
\be \label{LFD-nC2}
({\cal D}^{\alpha} f)(a):= 
\frac{1}{\Gamma(n-\alpha)} \lim_{x \to a}
\int^x_a \frac{dz}{(x-z)^{\alpha-n+1}} \frac{d^n f(z)}{dz^n},
\quad ( n-1 < \alpha \le n) .
\ee
}

{\bf Proof.}
This statement is direct corollary of Theorem 2.1 of \cite{KST} 
(see Eq. (2.4.15) of \cite{KST}), where
the left-sided Caputo fractional derivative is represented
in the form
\be \label{CFD-n2}
\, _a^{C}{\cal D}^{\alpha}_x f(x) :=
\frac{1}{\Gamma(n-\alpha)} 
\int^x_a \frac{dz}{(x-z)^{\alpha-n+1}} \frac{d^n f(z)}{dz^n},
\quad ( n-1 < \alpha \le n) .
\ee
\\

The Caputo fractional derivatives $\ _a^CD^{\alpha}_{x}$
can be defined for functions belonging to the space $AC^n[a,b]$
of absolutely continuous functions \cite{KST}.
If $f(x) \in AC^n[a,b]$, then the Caputo fractional derivatives
exist almost everywhere on $[a,b]$.

It can be directly derived that the left-sided Caputo 
fractional derivative of the power function $(x-a)^{\beta}$ is
\[ _a^C{\cal D}^{\alpha}_{x} (x-a)^{\beta}=
\frac{\Gamma(\beta+1)}{\Gamma(\beta-\alpha+1)} (x-a)^{\beta-\alpha} , \]
where $\beta>-1$ and $\alpha>0$.
In particular, if $\beta=0$ and $\alpha>0$,
then the Caputo fractional derivative of a constant $C$ 
is equal to zero:
\[  _a^C{\cal D}^{\alpha}_{x} C=0 . \]
For $k=0,1,2,...,n-1$, we have
\[ _a^C{\cal D}^{\alpha}_{x} (x-a)^k=0 . \]

It is important to emphasize that 
the Caputo fractional derivative 
is represented by equation (\ref{CFD-n2})
as a sequence of two operations, namely, first taking the ordinary derivative of integer order $n$ 
and then the integration of fractional order $n-\alpha$,
(see Eq. 2.4.17 of \cite{KST}):
\[ (\, _a^{C}{\cal D}^{\alpha}_x f)(x) =
\, _a^{RL}{\cal I}^{n-\alpha}_x f^{(n)}(x) .
\]
The Riemann-Liouville fractional derivative 
is represented by the inverse sequence of the same operations.


\section{Local fractional derivative is derivative of integer-order}

Let us prove the following statement for the local 
fractional derivatives of $n$-differentiable functions. \\

{\bf Proposition.} \\
{\it Local fractional derivative of order $\alpha$, 
where $n-1 < \alpha \le n$, $n\in \mathbb{N}$,
on a set of $n$-differentiable functions
is a derivative of order $n\in \mathbb{N}$ or zero operator.} \\

{\bf Proof.} \\
If $f(x)$ is a $n$-differentiable function in a neighborhood 
$U(x_0)$ of $x_0$
(i.e. $f(x)$ has a derivative of order $n \in \mathbb{N}$ 
in a neighborhood of $x_0$), then the function
can be represented by the Taylor series
\be \label{Dif-n}
f(x) = \sum^{n}_{k=0} 
\frac{f^{(k)}(x_0)}{k!} (x-x_0)^k + o((x-x_0)^{n}) .
\ee
If $f(x)$ is a differentiable function $(i.e. $n=1$)$, then
\be \label{Dif-1}
f(x) = f(x_0) +f^{(1)} (x-x_0) + o((x-x_0)) .
\ee
For $n$-differentiable functions, we can rewrite
the local fractional derivative  
of order $n-1 < \alpha \le n$ in the form
\[
({\cal D}^{\alpha} f)(x_0) = 
\lim_{x \to x_0} \, _{x_0}^{RL}{\cal D}^{\alpha}_x  
\left( f(x) - \sum^{n-1}_{k=0} \frac{f^{(k)}(x_0)}{k!} (x-x_0)^k \right) = 
\]
\[
=\lim_{x \to x_0}
\, _{x_0}^{RL}{\cal D}^{\alpha}_x  
\left( \frac{f^{(n)}(x_0)}{n!} (x-x_0)^n + o((x-x_0)^n)\right) = 
\]
\be \label{RLFD-n2}
= \frac{f^{(n)}(x_0)}{n!}  \lim_{x \to x_0}
\, _{x_0}^{RL}{\cal D}^{\alpha}_x  
\Bigl( (x-x_0)^n + o((x-x_0)^n) \Bigr) .
\ee
Using $\Gamma(n+1)=n!$ and the equation (see Property 2.1 of \cite{KST})
\be
 _{x_0}^{RL}{\cal D}^{\alpha}_x  (x-x_0)^{\beta} =
\frac{\Gamma(\beta+1)}{\Gamma(\beta+1-\alpha)}
(x-x_0)^{\beta-\alpha} , \quad (\beta >-1) ,
\ee
we obtain
\[
({\cal D}^{\alpha} f)(x_0) 
= \frac{f^{(n)}(x_0)}{n!}  \lim_{x \to x_0} 
\frac{\Gamma(n+1)}{\Gamma(n+1-\alpha)}
\Bigl( (x-x_0)^{n-\alpha} + o((x-x_0)^{n-\alpha}) \Bigr) = 
\]
\be \label{RLFD-n3}
= \frac{f^{(n)}(x_0)}{\Gamma(n+1-\alpha)}
\lim_{x \to x_0} \left( (x-x_0)^{n-\alpha} + o((x-x_0)^{n-\alpha})\right) .
\ee
As a result, we get
\be \label{RLFD-n4}
({\cal D}^{\alpha} f)(x_0) = 
\left\{
\begin{array}{cc}
0 & n-1<\alpha <n, \\
& \\
f^{(n)}(x_0) & \alpha =n ,
\end{array}
\right.
\ee
where $n-1 < \alpha \le n$ and $n \in \mathbb{N}$, 
and we use $\Gamma(n+1-n)=\Gamma (1)=1$. \\

The proved statement means that
the local fractional derivatives cannot be considered as
fractional derivatives on sets of finite (or infinite) numbers 
of times differentiable functions.

As a corollary of the proved statment we can say that 
generalizations of the Leibniz rule 
for local fractional derivatives ${\cal D}^{\alpha}$ 
of order $\alpha>0$ 
for $n$-differentiable functions $f(x)$ and $g(x)$
should give the relation 
\be \label{Req-2b}
{\cal D}^{\alpha} (f(x) \, g(x)) = \sum^n_{k=0} 
\frac{\Gamma(\alpha+1)}{\Gamma(\alpha-k+1) \, \Gamma(k+1)} 
({\cal D}^{\alpha-k} f(x)) \, ({\cal D}^{k} g(x)) , 
\quad (\alpha=n \in \mathbb{N}) 
\ee 
for integer orders $\alpha=n \in \mathbb{N}$ or
it should give
the equality $0=0$ for noninteger orders $\alpha$,
$n-1 < \alpha \le n$.

For $0<\alpha \le 1$, equation (\ref{RLFD-n4}) has the form
\be
({\cal D}^{\alpha} f)(x_0) = 
\left\{
\begin{array}{cc}
0 & 0<\alpha <1, \\
& \\
f^{\prime}(x_0) & \alpha =1 .
\end{array}
\right.
\ee

Therefore we get the statement that 
has been proved in \cite{CNSNS2013} as 
Theorem "No violation of the Leibniz rule. No fractional derivative":
{\it If the unviolated Leibniz rule
\be \label{uvLR} 
{\cal D}^{\alpha}  (f(x) \, g(x)) = 
({\cal D}^{\alpha} f(x)) \, g(x) + 
f(x) \, ({\cal D}^{\alpha}  g(x)) , 
\ee
holds for the linear operator ${\cal D}^{\alpha}$ 
and the $n$-differentiable functions $f(x)$ and $g(x)$,
then this operator is the derivative 
of integer (first) order, that can 
be represented as ${\cal D}^{\alpha} = a(x) \, D^1$,
where $a(x)$ are functions on $\mathbb{R}$.} \\

In paper \cite{CNSNS2013} it was proved
that violation of relation (\ref{uvLR}) 
is a characteristic property of 
all types of fractional-order derivatives with $\alpha > 0$
(and it is obvious for
derivatives of integer-orders $n \in \mathbb{N}$ 
greater than one).

A correct form of the Leibniz rule for local
fractional derivatives should be obtained from
the fractional generalization
of the Leibniz rule for the Riemann-Liouville derivatives
that has the form of the infinite series.
The symmetrized expression of this rule is
\be \label{GLR}
_a^{RL}{\cal D}^{\alpha}_x (fg)(x) =
\sum^{\infty}_{k=0} \binom{k}{\alpha} \Bigl(
( \, _a^{RL}{\cal D}^{\alpha-k}_x f)(x) \, (D^{k} g(x)) +
( \, _a^{RL}{\cal D}^{\alpha-k}_x g)(x) \, (D^{k} f(x)) 
\Bigr),
\ee
where the $f(x)$ and $g(x)$ are analytic functions 
on the interval $[a,b]$ (Theorem 15.1 of \cite{SKM}), 
$D^k f(x)= f^{(k)}(x)$, and
\[ 
\binom{k}{\alpha} = \frac{\Gamma(\alpha+1)}{2\, 
\Gamma(\alpha-k+1) \, \Gamma(k+1) }  . 
\]
It should be noted that
the sum of (\ref{GLR}) is infinite and it contains 
integrals of fractional order for the values $k>[\alpha]+1$.
Equation (\ref{GLR}) has been proposed  by Liouville
\cite{Liouv} in 1832.
Note that extensions of the Leibniz rule to fractional
derivatives are suggested by Osler 
\cite{Osler1970,Osler1971,Osler1972,Osler1973}.

The left-sided Riemann-Liouville fractional derivative
of order $\alpha$ ($n-1< \alpha\le n$) 
of a function $f(x)$ (see Eqs. (2.4.6) of \cite{KST})
can be presented via the Caputo fractional derivative by
\be \label{RLFD-CFD}
(\, _a^{RL}{\cal D}^{\alpha}_x f)(x) = 
(\, _a^{C}{\cal D}^{\alpha}_x  f)(x)
+ \sum^{n-1}_{k=0} \frac{f^{(k)}(a)}{k!} (x-a)^{k-\alpha} , 
\quad ( n-1< \alpha \le n) .
\ee
Substitution of equation (\ref{RLFD-CFD})
for $(\, _a^{RL}{\cal D}^{\alpha} fg)(x)$, 
$(\, _a^{RL}{\cal D}^{\alpha} f)(x)$ and
$(\, _a^{RL}{\cal D}^{\alpha} g)(x)$
into equations (\ref{GLR}) gives
a generalized Leibniz rule for the Caputo fractional derivatives.
The extension of the Leibniz rule for
the Caputo fractional derivatives
also has been suggested by Theorem 3.17 of the Diethelm's 
book \cite{Diethelm}.


Let us give some remarks about Leibniz rule (\ref{uvLR}) 
for local fractional derivatives on a set of 
non-differentiable functions.
It should be noted that statement \cite{BAC,BG}
that the Leibniz rule in the unviolated form (\ref{uvLR}) 
holds for non-differentiable functions is incorrect 
\cite{PLA2015}. 
Let us give some explanations below.

(1) It should be noted that the Leibniz rule (\ref{uvLR}) 
cannot be used for non-differentiable functions 
that are not fractional-differentiable
since the derivatives 
${\cal D}^{\alpha} f(x)$,  
${\cal D}^{\alpha}  g(x)$
and  ${\cal D}^{\alpha}  (f(x) \, g(x))$
should exist.

(2) It is easy to see that nowhere in the proofs 
proposed in \cite{BAC,BG}, the requirement 
that the functions  $f(x)$ and $g(x)$ are not 
classically differentiable is not used.
Therefore, using the same proofs, 
we can get the statement that the Leibniz rule (\ref{uvLR}) 
holds for fractional-differentiable functions $f(x)$, $g(x)$ 
without the useless assumption that these functions 
are not classically differentiable. 

(3) The Caputo fractional derivatives $\ _a^CD^{\alpha}_{x}$ 
can be defined for functions belonging to the space $AC^n[a,b]$
of absolutely continuous functions \cite{KST}.
If $f(x) \in AC^n[a,b]$, then the Caputo fractional derivatives
exist almost everywhere on $[a,b]$.
If $f(x) \in AC^n[a,b]$, then this functions
is $n$-differentiable almost everywhere on $[a,b]$.
The Leibniz rule (\ref{uvLR}) holds on
a set of differentiable functions only for $\alpha=1$.

As a result, we get that the Leibniz rule (\ref{uvLR})
with local fractional derivative 
for $f(x), g(x) \in AC^n[a,b]$ holds only 
for $\alpha=1$ almost everywhere on $[a,b]$.
In the countable number of points the Leibniz rule 
is not performed since the derivatives 
${\cal D}^{\alpha} f(x)$,  
${\cal D}^{\alpha}  g(x)$
and  ${\cal D}^{\alpha}  (f(x) \, g(x))$
does not exist.
Note that the Leibniz rule 
for fractional derivatives of orders $\alpha =n>1$ 
($n\in \mathbb{N}$) cannot have the form (\ref{uvLR}).
It should be get the rule in the form (\ref{Req-2b}).

\section{Conclusion}

We prove that 
local fractional derivatives of differentiable functions 
are integer-order derivatives or zero operator.
We show that definition of local fractional derivative 
of non-integer order $\alpha$ ($n-1 <\alpha \le n$)
can be represented by limit of 
the Caputo fractional derivative of order $\alpha$.
The Caputo derivative of fractional order $\alpha$
of a function $f(x)$ is defined 
as an action of derivative of integer order $n$
and then the integration of non-integer order $n-\alpha$.
Therefore the function $f(x)$ 
should be classically differentiable $n$ times, and
derivative $f^{(n)}(x)$ of integer order $n$ should exist.
As a result, the local fractional derivative cannot 
be considered as the best tool for description of
nowhere differentiable functions and fractal objects
despite the motivation to consider these derivatives.
We also argued that unviolated Leibniz rule cannot hold 
for local fractional derivatives of orders $\alpha \ne 1$.

A problem of strict formulation of
"fractal" analogs of derivatives 
and integrals by extension of the classical and 
fractional calculus to
nowhere differentiable functions and fractal sets
is very difficult.
We can assume that the most promising approaches
to fractional calculus on fractals can be based
on methods of "non-integer dimensional spaces"
(see \cite{CNSNS2015,JMP2014} and references therein)
or "analysis on fractals" 
\cite{Kugami,Strichartz-1,Harrison,Kumagai}.



\end{document}